\documentclass[a4paper,10pt, final]{article}

\usepackage[notref,notcite]{showkeys}
\usepackage{amsthm,amsfonts,amsmath}
\usepackage{color}
\usepackage{graphicx,psfrag}
\usepackage{hyperref}
\usepackage{a4wide}
\usepackage{enumerate,float}

\usepackage{pstricks}
\usepackage{pstricks-add}
\usepackage{epic,eepic}
\usepackage{pst-math,pst-plot}
\usepackage{pst-grad, pst-text, pst-node, pst-tree}

\newtheorem{proposition}{Proposition}[section]
\newtheorem{theorem}[proposition]{Theorem}

\newtheorem{definition}[proposition]{Definition}

\newenvironment{proofof}[1]{\smallskip\noindent{\textbf{Proof~of~#1.}}%
  \hspace{1pt}}{\hspace{-5pt}{\nobreak\quad\nobreak\hfill\nobreak%
    $\square$\vspace{2pt}\par}\smallskip\goodbreak}

\numberwithin{equation}{section}
\allowdisplaybreaks[4]

\renewcommand{\phi}{\varphi}
\renewcommand{\epsilon}{\varepsilon}
\renewcommand{\theta}{\vartheta}
\renewcommand{\L}[1]{\mathbf{L^#1}}
\newcommand{\Lloc}[1]{\mathbf{L^{#1}_{loc}}}
\newcommand{\C}[1]{\mathbf{C^{#1}}}
\newcommand{\Cc}[1]{\mathbf{C_c^{#1}}}

\newcommand{\W}[2]{\mathbf{W^{#1,#2}}}
\newcommand{\BV}{\mathbf{BV}}
\newcommand{\modulo}[1]{{\left|#1\right|}}
\newcommand{\norma}[1]{{\left\|#1\right\|}}
\newcommand{\reali}{{\mathbb{R}}}

\newcommand{\tv}{\mathop\mathrm{TV}}

\newcommand{\Id}{\mathop{\mathbf{Id}}}

\renewcommand{\d}[1]{\mathinner{\mathrm{d}{#1}}}

\begin{document}

\title{The Compressible to Incompressible Limit \\ of 1D Euler
  Equations: the Non Smooth Case}

\author{Rinaldo M. Colombo$^1$ \and Graziano Guerra$^2$ \and Veronika
  Schleper$^3$}

\footnotetext[1]{INDAM Unit, University of Brescia,
  Italy. \texttt{rinaldo.colombo@unibs.it}}

\footnotetext[2]{Department of Mathematics and Applications,
  University of Milano - Bicocca,
  Italy. \texttt{graziano.guerra@unimib.it}}

\footnotetext[3]{Institute for Applied Analysis and Numerical
  Simulations, University of Stuttgart,
  Germany. \texttt{veronika.schleper@mathematik.uni-stuttgart.de}}

\maketitle

\begin{abstract}
  \noindent We prove a rigorous convergence result for the
  compressible to incompressible limit of weak entropy solutions to
  the isothermal 1D Euler equations.

  \medskip


  \medskip

\end{abstract}

\section{Introduction}
\label{sec:Intro}

The compressible to incompressible limit is widely studied in the
literature, see for instance the well known
results~\cite{KlainermanMajda1981, KlainermanMajda1982,
  MetivierSchochet2001, Schochet1986, Secchi2000}, the more
recent~\cite{JiangYong}, the review~\cite{Schochet2005} and the
references therein. The usual setting considers regular solutions,
whose existence is proved only for a finite time, to the compressible
equations in $2$ or $3$ space dimensions. As the Mach number vanishes,
these solutions are proved to converge to the solutions of the
incompressible Euler equations, while only a weak* convergence holds
for the pressure.

Here, we recover the same convergence results, in a 1D setting, but
within the framework of \emph{weak entropy solutions} proved to exist
for \emph{all times}. In particular, $\L1$-convergence is proved for
all positive times and the space regularity of solutions is $\L1 \cap
\BV$.

In the present 1D setting, a compressible to incompressible limit over
all the real line is of no interest. One may then consider two
compressible immiscible fluids, say a gas and a liquid, letting the
liquid tend to become incompressible. The limiting procedure basically
yields a boundary value problem for the gas, while the boundary
between the fluids turns into a solid wall. Indeed, an incompressible
fluid in 1D is essentially a solid.

Therefore, below we consider a droplet of a compressible inviscid
fluid that fills the segment $[a,b]$ and is surrounded by another
compressible inviscid fluid filling the rest of the real line. The
fluids are assumed to be immiscible. For simplicity, we refer to the
fluid forming the droplet as to a \emph{liquid}, while its complement
is labeled as \emph{gas}. In the isentropic (or isothermal)
approximation, this whole system can be described by
\begin{equation}
  \label{eq:Model}
  \left\{
    \begin{array}{lr@{\,}c@{\,}l@{\qquad}l@{}}
      \left\{
        \begin{array}{@{}l}
          \partial_t \rho_g + \partial_x (\rho_g \, v_g) = 0
          \\
          \partial_t (\rho_g \, v_g)
          +
          \partial_x \left( \rho_g \, {v_g}^2 + p_g(\rho_g) \right)
          =
          0
        \end{array}
      \right.
      &
      x & \in & \reali\setminus [a (t), b (t)]
      &
      \mbox{gas}
      \\[15pt]
      \left\{
        \begin{array}{@{}l}
          \partial_t \rho_l + \partial_x (\rho_l \, v_l) = 0
          \\
          \partial_t (\rho_l \, v_l)
          +
          \partial_x \left( \rho_l \, {v_l}^2 + p_l(\rho_l) \right)
          =
          0
        \end{array}
      \right.
      &
      x & \in & [a (t), b(t)]
      &
      \mbox{liquid}
      \\[15pt]
      \left\{
        \begin{array}{rcl}
          \dot a (t) & = & v_g\left(t, a (t)-\right)
          \\
          \dot b (t) & = & v_g\left(t, b (t)+\right)
        \end{array}
      \right.
      &
      & &
      &
      \begin{array}{@{}l@{}}
        \mbox{kinetic}
        \\
        \mbox{relations}
        \\
        \mbox{(immiscibility)}
      \end{array}
      \\[15pt]
      \left\{
        \begin{array}{@{}rcl@{}}
          v_g \left(t, a (t)-\right)
          & = &
          v_l \left(t, a (t)+\right)
          \\
          v_g \left(t, b (t)+\right)
          & = &
          v_l \left(t, b (t)-\right)
        \end{array}
      \right.
      &
      & &
      &
      \begin{array}{@{}l@{}}
        \mbox{mass}\\\mbox{conservation}
      \end{array}
      \\[15pt]
      \left\{
        \begin{array}{rcl}
          p_g\left(\rho_g(t, a (t)-) \right)
          & = &
          p_l \left(\rho_l(t, a (t)+)\right)
          \\
          p_g\left(\rho_g(t, b (t)+) \right)
          & = &
          p_l \left(\rho_l(t, b (t)-)\right)
        \end{array}
      \right.
      &
      & &
      &
      \begin{array}{@{}l@{}}
        \mbox{momentum}\\\mbox{conservation.}
      \end{array}
    \end{array}
  \right.
\end{equation}
Here, $\rho_l$ is the density of liquid in the droplet, $v_l$ is its
speed and $p_l = p_l (\rho_l)$ is the pressure law, while $\rho_g$,
$v_g$ and $p_g$ denote the analogous quantities for the gas. The above
conditions not only ensure the conservation of mass, but also prevent
any exchange of matter between the two fluids. In particular, they do
not mix but, clearly, there is exchange of information between the two
fluids, thanks to the global conservation of momentum. Note also that
these conditions obviously ensure the energy conservation at the
interfaces, while shocks lead to the dissipation of energy in the
interior of the 2 fluids. Indeed, as is well know~\cite{DafermosBook},
energy plays here the role of the mathematical entropy. The above
kinetic relations provide a link between the fluid interfaces and the
fluid speeds.


Passing to the incompressible limit in the liquid phase, we expect to
obtain the system
\begin{equation}
  \label{eq:Incomp}
  \left\{
    \begin{array}{lr@{\,}c@{\,}l@{\qquad\qquad}l@{}}
      \left\{
        \begin{array}{@{}l}
          \partial_t \rho_g + \partial_x (\rho_g \, v_g) = 0
          \\
          \partial_t (\rho_g \, v_g)
          +
          \partial_x \left( \rho_g \, {v_g}^2 + p_g(\rho_g) \right)
          =
          0
        \end{array}
      \right.
      &
      x & \in & \reali\setminus [a (t), b (t)]
      &
      \mbox{gas;}
      \\[15pt]
      \left\{
        \begin{array}{@{}l}
          \rho_l = \bar\rho_l
          \\
          \dot v_l
          =
          \frac{p_g\left(t, a (t-)\right)- p_g\left(t, b (t+)\right)}{\left(b (t)-a (t)\right) \, \bar\rho_l}
        \end{array}
      \right.
      &
      x & \in & [a (t), b(t)]
      &
      \mbox{liquid;}
      \\[15pt]
      \left\{
        \begin{array}{@{}rcl@{}}
          v_g \left(t, a (t)-\right)
          & = &
          \dot a (t)
          \\
          v_g \left(t, b (t)+\right)
          & = &
          \dot b (t)
          \\
          v_l (t)\; = \; \dot a (t) & = & \dot b (t)
        \end{array}
      \right.
      &
      & &
      &
      \begin{array}{@{}l}
        \mbox{mass and}
        \\
        \mbox{momentum}
        \\
        \mbox{conservation,}
      \end{array}
    \end{array}
  \right.
\end{equation}
consisting of a conservation law describing the compressible gas
coupled with an ordinary differential equation for the incompressible
droplet. We recall that~\eqref{eq:Incomp} is known to be well posed,
see~\cite[Proposition~3.1]{BorscheColomboGaravello2}.

Since we assume the two phases immiscible, a natural choice is to pass
to Lagrangian coordinates, so that the interfaces between the two
phases become stationary.

In these coordinates, introducing the total mass $m$ of the liquid, the
above system~\eqref{eq:Model} reads:
\begin{equation}
  \label{eq:ModelL}
  \left\{
    \begin{array}{lr@{\,}c@{\,}l@{\qquad}l@{}}
      \left\{
        \begin{array}{@{}l@{}}
          \partial_t \tau_g - \partial_z v_g = 0
          \\
          \partial_t v_g + \partial_z p_g (\tau_g) =0
        \end{array}
      \right.
      &
      z & \in & \reali\setminus [0, m]
      &
      \mbox{gas}
      \\[15pt]
      \left\{
        \begin{array}{@{}l@{}}
          \partial_t \tau_l - \partial_z v_l = 0
          \\
          \partial_t v_l + \partial_z p_l (\tau_l) =0
        \end{array}
      \right.
      &
      z & \in & [0, m]
      &
      \mbox{liquid}
      \\[15pt]
      \left\{
        \begin{array}{@{}rcl@{}}
          v_g \left(t, 0-\right)
          & = &
          v_l \left(t, 0+\right)
          \\
          v_g \left(t, m+\right)
          & = &
          v_l \left(t, m-\right)
        \end{array}
      \right.
      &
      & &
      &
      \begin{array}{@{}l@{}}
        \mbox{mass}\\\mbox{conservation}
      \end{array}
      \\[15pt]
      \left\{
        \begin{array}{rcl}
          p_g\left(\tau_g(t, 0-) \right)
          & = &
          p_l \left(\tau_l(t, 0+)\right)
          \\
          p_g\left(\tau_g(t, m+) \right)
          & = &
          p_l \left(\tau_l(t, m-)\right)
        \end{array}
      \right.
      &
      & &
      &
      \begin{array}{@{}l}
        \mbox{momentum}
        \\
        \mbox{conservation,}
      \end{array}
    \end{array}
  \right.
\end{equation}
while its incompressible limit in the liquid phase,
i.e.,~\eqref{eq:Incomp} in Lagrangian coordinates, becomes
\begin{equation}
  \label{eq:IncompL}
  \left\{
    \begin{array}{lr@{\,}c@{\,}l@{\qquad\qquad}l@{}}
      \left\{
        \begin{array}{@{}l}
          \partial_t \tau_g - \partial_z v_g = 0
          \\
          \partial_t v_g
          +
          \partial_z  p_g(\tau_g)
          =
          0
        \end{array}
      \right.
      &
      z & \in & \reali\setminus [0, m]
      &
      \mbox{gas;}
      \\[20pt]
      \dot v_l
      =
      \frac{p_g(t, 0-) - p_g(t, m+)}{m}
      & & & &
      \mbox{liquid;}
      \\[10pt]
      \left\{
        \begin{array}{@{}rcl@{}}
          v_g \left(t, 0-\right)
          & = &
          v_l (t)
          \\
          v_g \left(t, m+\right)
          & = &
          v_l (t)
        \end{array}
      \right.
      &
      & &
      &
      \begin{array}{@{}l}
        \mbox{mass and momentum}
        \\
        \mbox{conservation.}
      \end{array}
    \end{array}
  \right.
\end{equation}
Problems~\eqref{eq:Model} and~\eqref{eq:ModelL},
respectively~\eqref{eq:Incomp} and~\eqref{eq:IncompL}, are well posed
in $\L1$ globally in time, see~\cite[Theorem~2.5]{ColomboSchleper},
respectively~\cite[Theorem~3.6]{BorscheColomboGaravello3}.

The main result of this paper states the rigorous convergence
of~\eqref{eq:Model} to~\eqref{eq:Incomp} or, equivalently,
of~\eqref{eq:ModelL} to~\eqref{eq:IncompL}. In this very singular
limit, the sound speed $\bar c$ in the liquid phase tends to
$+\infty$; the density converges to a fixed reference value
$\bar\rho$; the graph of the pressure law becomes vertical and the
eigenvectors of the Jacobian of the flow tend to
coalesce. In~\cite{KlainermanMajda1981, KlainermanMajda1982,
  MetivierSchochet2001, Schochet1986, Schochet2005}, this behavior is
reproduced requiring that $\norma{\nabla p} \to +\infty$. Here,
similarly, we fix the constant initial specific volume $\bar\tau_l$ of
the liquid and locally describe the liquid phase in Lagrangian
coordinates in a neighborhood of $\bar\tau$ through the pressure law
\begin{equation}
  \label{eq:linearPressure}
  p^\eta (\tau) = \bar p - \eta^2 (\tau - \bar\tau_l)
  \quad \mbox{ so that } \quad
  \tau^\eta (p) = \bar\tau_l - \frac{1}{\eta^2} \, (p-\bar p)
  \qquad \mbox{ where } \qquad
  \eta = \frac{\bar c}{\bar\tau_l} \,,
\end{equation}
where, $\bar c$ is the Eulerian sound speed and $\eta$ is the
Lagrangian sound speed.

\medskip

Informally, we present here the main result of the present paper. The
rigorous statement is in Section~\ref{sec:WP}. More precisely,
points~1. and 2.~in Theorem~\ref{thm:one} are consequences of
Theorem~\ref{thm:WFT}, while point 3.~follows from
Theorem~\ref{thm:Main}.

\begin{theorem}
  \label{thm:one}
  Fix in the gas phase a pressure law satisfying~\textbf{(p)} and the
  pressure law~\eqref{eq:linearPressure} in the liquid phase.
  \begin{enumerate}
  \item There exists a nontrivial set of initial data $\mathcal{D}$
    such that the Cauchy problem consisting of~\eqref{eq:Model} is
    well posed for all initial data in $\mathcal{D}$ and for all $\bar
    c$ sufficiently large.
  \item The initial data in $\mathcal{D}$ in the liquid phase are
    incompressible, i.e., they have a constant density $\bar\rho_l$
    and speed $\bar v_l$.
  \item For any initial data in $\mathcal{D}$, the solution to the
    Cauchy problem~\eqref{eq:Model} converges as $\bar c\to +\infty$
    and the limit solves~\eqref{eq:Incomp}. In particular, in the
    liquid phase:
    \begin{enumerate}[(i)]
    \item the density converges in $\Lloc1$ to $\bar \rho_l$, which is
      constant in space and time;
    \item the velocity converges in $\Lloc1$ to a function $v_l =
      v_l(t)$ which is constant in space and whose time evolution is
      specified in~\eqref{eq:Incomp};
    \item the pressure weak* converges in $\L\infty$ to the linear
      interpolation of the values of the gas pressure at the liquid
      boundaries, for almost every time $t\geq 0$.
    \end{enumerate}
  \end{enumerate}
\end{theorem}

\noindent Remark that the convergences stated above extend the
classical result~\cite{KlainermanMajda1981} to the case of \emph{non
  smooth}, albeit 1D, solutions.

The present setting comprises also the case of the compressible to
incompressible limit on the free boundary value problem in which only
the liquid is present and it is constrained in the moving strip $[a
(t), b (t)]$, with prescribed pressure along the boundaries.

\medskip

As usual in the study of non smooth solutions to 1D systems of
conservation laws, see~\cite{BressanLectureNotes, DafermosBook}, we
exploit the wave front tracking techniques, so that the desired
estimates follow from suitable bounds on the approximate solutions. A
first key analytical difficulty in obtaining the present result lies
in the need for bounds in the total variation of solutions that are
uniform when the sound speed tends to $+\infty$. A very careful choice
of the parametrization allows to control the interactions of waves
against the phase boundaries.

Within the liquid phase, waves may well bounce back and forth with a
diverging speed. This phenomenon requires an \emph{ad hoc} procedure
to bound the total number of interaction points in the
$\epsilon$-approximations. In turn, we are able to use the wave front
tracking algorithm in~\cite{BressanLectureNotes}, but without the need
of nonphysical waves. A further consequence of these possible bounces
is that the total variation of the pressure along lines $x = \bar x$
grows unboundedly, which is why only a weak* convergence of the
pressure is possible. Nevertheless, we recover the Newton law for the
incompressible liquid thanks to the bounds on the total variation and
to the conservative form of~\eqref{eq:ModelL} and~\eqref{eq:IncompL}.

\medskip

The next Section presents the rigorous setting for~\eqref{eq:ModelL}
and~\eqref{eq:IncompL}. First, the well posedness of~\eqref{eq:ModelL}
proved in~\cite[Theorem~2.5]{ColomboSchleper} is improved to obtain
uniform estimates that allow to pass to the incompressible
limit. Then, Theorem~\ref{thm:Main} presents our main
result. Section~\ref{sec:TD} is devoted to the analytical
proofs. There, we highlight the key modifications to the standard
procedure~\cite{BressanLectureNotes, DafermosBook}, without repeating
the now classical wave front tracking constructions.

\section{The Compressible \texorpdfstring{$\to$}{to} Incompressible Limit}
\label{sec:WP}

Throughout, we identify the state $u$ of the fluid by $(\tau,v)$ or
$(p,v)$, once an invertible pressure law $p = p (\tau)$ is
assigned. Given a function $u \colon \reali \to \reali^+ \times
\reali$, we denote
\begin{displaymath}
   (\tau_g, v_g) = u_g = u \, \chi_{\strut \reali \setminus [0, m]}
  \qquad \mbox{ and } \qquad
  (\tau_l, v_l) = u_l = u \, \chi_{\strut[0, m]} \,.
\end{displaymath}
In Lagrangian coordinates, a standard assumption on the pressure law
is
\begin{description}
\item[(p)] $p \in \C4 (\reali^+; \reali^+)$ is such that $p' (\tau)
  <0$ and $p'' (\tau) > 0$ for all $\tau>0$.
\end{description}

\noindent Recall the definition of solution or~\eqref{eq:ModelL}:

\begin{definition}
  \label{def:solution}
  Fix a gas state $\bar u_g$ and a liquid state $\bar u_l$. Denote
  $\bar u = \bar u_g + (\bar u_l - \bar u_g)\, \chi_{\strut[0,m]}$.
  Let $T>0$ be fixed. By solution to~\eqref{eq:ModelL} we mean a map
  \begin{displaymath}
    u
    \in \C0 \left([0,T]; \bar u + (\L1 \cap \BV) (\reali; \reali^+
      \times \reali)\right)
  \end{displaymath}
  such that:
  \begin{enumerate}
  \item it is a weak entropy solution to $\left\{
      \begin{array}{@{}l@{}}
        \partial_t \tau_g - \partial_z v_g = 0
        \\
        \partial_t v_g + \partial_z p_g (\tau_g) =0
      \end{array}
    \right.$ in $[0,T] \times (\reali \setminus [0,m])$;
  \item it is a weak entropy solution to $\left\{
      \begin{array}{@{}l@{}}
        \partial_t \tau_l - \partial_z v_l = 0
        \\
        \partial_t v_l + \partial_z p_l (\tau_l) =0
      \end{array}
    \right.$ in $[0,T] \times [0, m]$;
  \item for a.e.~$t \in [0,T]$, the conditions at the junction
    $\left\{
      \begin{array}{@{}rcl@{}}
        p_g\left(\tau_g(t, 0-) \right)
        & = &
        p_l \left(\tau_l(t, 0+)\right)
        \\
        v_g \left(t, 0-\right)
        & = &
        v_l \left(t, 0+\right)
      \end{array}
    \right.$ and $\left\{
      \begin{array}{@{}rcl@{}}
        p_g\left(\tau_g(t, m+) \right)
        & = &
        p_l \left(\tau_l(t, m-)\right)
        \\
        v_g \left(t, m-\right)
        & = &
        v_l \left(t, m+\right)
      \end{array}
    \right.$ are satisfied.
  \end{enumerate}
\end{definition}

\noindent Above, the continuity in $\C0 \left([0,T]; \bar u + (\L1
  \cap \BV) (\reali; \reali^+ \times \reali)\right)$ is understood
with respect to the $\L1$ norm.

In the case of the mixed model~\eqref{eq:IncompL}, we
adapt~\cite[Definition~3.1]{BorscheColomboGaravello3} to the present
situation.

\begin{definition}
  \label{def:solInc}
  Fix a gas state $\bar u_g$. Let $T>0$ be fixed. By solution
  to~\eqref{eq:IncompL} we mean a pair of functions
  \begin{equation}
    \label{eq:regularity}
    (u, v_l) \in
    \C0\left(
      [0,T];
      \bar u_g
      +
      (\L1 \cap \BV) (\reali \setminus [0, m]; \reali^+ \times \reali)
    \right)
    \times
    \W{1}{\infty} ([0,T]; \reali)
  \end{equation}
  such that:
  \begin{enumerate}
  \item $u$ is a weak entropy solution to $\left\{
      \begin{array}{@{}l@{}}
        \partial_t \tau_g - \partial_z v_g = 0
        \\
        \partial_t v_g + \partial_z p_g (\tau_g) =0
      \end{array}
    \right.$ in $[0,T] \times (\reali \setminus [0,m])$;
  \item $v_l$ is a solution to $\dot v_l = \frac{1}{m}\left( p_g
      \left(\tau_g(s, 0-)\right) - p_g \left(\tau_g(s, m+)\right)
    \right)$ on $[0,T]$;
  \item the conditions at the interface $\left\{
      \begin{array}{@{}rcl@{}}
        v_g \left(t, 0-\right)
        & = &
        v_l (t)
        \\
        v_g \left(t, m+\right)
        & = &
        v_l (t)
      \end{array}
    \right.$ are satisfied for a.e.~$t \in [0,T]$.
  \end{enumerate}
\end{definition}

\noindent The definitions of solutions to the Cauchy problems
for~\eqref{eq:ModelL} and~\eqref{eq:IncompL} are an immediate
adaptation of the definitions above.

Aiming at the rigorous compressible $\to$ incompressible limit, we
need new estimates on the compressible problem~\eqref{eq:ModelL},
improving analogous results in~\cite[Theorem~2.5]{ColomboSchleper}.

\begin{theorem}
  \label{thm:WFT}
  Fix a state $\bar u_g$ in the gas phase and a state $\bar u_l$ in
  the liquid phase. Let $p_g$ satisfy~\textbf{(p)} and $p^\eta$ be as
  in~\eqref{eq:linearPressure}. Then, there exist positive $\Delta,
  \delta_g, \eta_*, L$ such that for all $\eta > \eta_*$ and for
  suitable positive $\delta_l^\eta$, $L^\eta$,
  problem~\eqref{eq:ModelL}--~\eqref{eq:IncompL} generates a semigroup
  \begin{displaymath}
    S^\eta \colon \reali^+ \times \mathcal{D}^\eta \to \mathcal{D}^\eta
  \end{displaymath}
  with the following properties:
  \begin{enumerate}
    \setlength{\parskip}{1pt} \setlength{\itemsep}{1pt}
  \item \label{it:thm1} $\mathcal{D}^\eta \supseteq \left\{u \in \bar
      u + (\L1\cap\BV) (\reali; \reali^+ \times \reali)\colon \tv
      (u_g) < \delta_g \mbox{ and } \tv (u_l) < \delta_l^\eta
    \right\}$;
  \item \label{it:thm2} $S^\eta$ is a semigroup: $S^\eta_0 = \Id$ and
    $S^\eta_{t_1} \circ S^\eta_{t_2} = S^\eta_{t_1+t_2}$;
  \item \label{it:thm:Lip} $S^\eta$ is Lipschitz in $u$: for any $u^1,
    u^2$ in $\mathcal{D}^\eta$ and for all $t \in \reali^+$
    \begin{displaymath}
      \norma{S^\eta_t u^1 - S^\eta_t u^2}_{\L1}
      \leq
      L^\eta \cdot \norma{u^1 - u^2}_{\L1} \,;
    \end{displaymath}
  \item \label{it:thm:LipT} $S^\eta$ is Lipschitz in $t$: for any $u$
    in $\mathcal{D}^\eta$ and for all $t_1,t_2 \in \reali^+$, setting
    $(\tau^\eta,v^\eta) (t) = S^\eta_t u$
    \begin{eqnarray*}
      \norma{\tau^\eta_g (t_1) - \tau^\eta_g (t_2)}_{\L1 (\reali\setminus[0,m])}
      & \leq &
      L \, \modulo{t_1-t_2}
      \\
      \norma{\tau^\eta_l (t_1) - \tau^\eta_l (t_2)}_{\L1 ([0,m])}
      & \leq &
      \frac{1}{\eta} \, L \, \modulo{t_1-t_2}
      \\
      \norma{v^\eta (t_1) - v^\eta (t_2)}_{\L1 (\reali)}
      & \leq &
      L \, \modulo{t_1-t_2} \,;
    \end{eqnarray*}
  \item \label{it:thm4} if $u \in \mathcal{D}^\eta$ is piecewise
    constant then, for $t$ small, the map $S^\eta_t u$ locally
    coincides with the standard Lax solutions to the Riemann problems
    for~\eqref{eq:ModelL} at the points of jump, at $x=0$ and at
    $x=m$;
  \item \label{it:thm5} for all $u \in \mathcal{D}^\eta$, the map $t
    \to S^\eta_t u$ is an entropy solution to~\eqref{eq:ModelL} with
    initial datum $u$ in the sense of Definition~\ref{def:solution};
  \item \label{it:thm:7} if $u \in \mathcal{D}^\eta$, then for all $t
    \geq 0$ and $\eta \geq \eta_*$, calling $\left(\tau^\eta (t),
      v^\eta (t)\right) = S^\eta_tu$
    \begin{equation}
      \label{eq:new}
      \begin{array}{rcl}
        \tv \left(p^\eta \left(\tau_l^\eta (t, \cdot)\right)\right)
        +
        \tv \left(p_g \left(\tau_g^\eta (t,\cdot)\right)\right)
        & < &
        \Delta
        \\[10pt]
        \eta \tv \left(v_l^\eta (t, \cdot)\right)
        +
        \eta^{2} \tv \left(\tau_l^\eta (t,\cdot)\right)
        & < &
        \Delta \,;
      \end{array}
    \end{equation}
  \item \label{it:thm:TVt} if $u \in \mathcal{D}^\eta$, then for all
    $\eta \geq \eta_*$, calling $\left(\tau^\eta (t), v^\eta
      (t)\right) = S^\eta_tu$
    \begin{equation}
      \label{eq:1}
      \begin{array}{rcl@{\qquad}r@{\,}c@{\,}l}
        \tv\left(\tau_g^\eta (\cdot, x)\right)
        +
        \tv\left(v_g^\eta (\cdot, x)\right)
        +
        \tv\left(p_g\left(\tau_g^\eta (\cdot, x)\right)\right)
        & \leq &
        \Delta
        &
        x & \in & \reali \setminus [0,m] \,;
      \end{array}
    \end{equation}
  \item \label{it:thm:last} for a.e.~$x_1,x_2$ with either $x_1,x_2 <
    0$ or $x_1,x_2>m$
    \begin{equation}
      \label{eq:uno}
      \int_0^t \norma{u_g^\eta (s, x_2) - u_g^\eta (s,x_1)} \d{s}
      \leq
      L \, \modulo{x_2-x_1} \,.
    \end{equation}
  \end{enumerate}
\end{theorem}

\noindent The points~\ref{it:thm2}., \ref{it:thm:Lip}.,
\ref{it:thm4}.~and~\ref{it:thm5}.~follow
from~\cite[Theorem~2.5]{ColomboSchleper}. To obtain the \emph{a
  priori} estimates in points~\ref{it:thm1}.,
\ref{it:thm:7}.~and~\ref{it:thm:TVt}., we have to substantially improve
the wave front tracking construction in~\cite{ColomboSchleper},
devising and exploiting a different parametrization of the Lax
curves. These estimates allow to obtain the key Lipschitz type
estimates at point~\ref{it:thm:LipT}.~and~\ref{it:thm:last}., where the dependence of the Lipschitz constants on $\eta$, where present, is explicit.

Remark that
point~\ref{it:thm:TVt}.~above may not hold in the liquid phase, since
the total variation of the pressure therein may well blow up. The
proof of Theorem~\ref{thm:WFT} is deferred to Section~\ref{sec:TD}.

The bounds above have a key role in proving our main result
below.

\begin{theorem}
  \label{thm:Main}
  Fix a state $\bar u_g$ in the gas phase and a state $\bar u_l$ in
  the liquid phase. Let $p$ satisfy~\textbf{(p)} and $p^\eta$ be as
  in~\eqref{eq:linearPressure}.  Then, with reference to the semigroup
  $S^\eta$ and its domain $\mathcal{D}^\eta$ constructed in
  Theorem~\ref{thm:WFT}, there exist positive $\delta_*$ and $\eta_*$
  such that for all $u_g^o \in \bar u_g^o + \L1(\reali\setminus [0,
  m]; \reali^+ \times \reali)$ with $\tv (u_g^o) < \delta_*$:
  \begin{enumerate}
  \item the function $u^o = u_g^o \, \chi_{\strut
      \reali\setminus[0,m]} + \bar u_l \, \chi_{\strut [0,m]}$ is in
    $\mathcal{D}^\eta$ for all $\eta > \eta_*$;
  \item the limit $(\tau_*, v_*) = \lim_{\eta \to +\infty} S^\eta_t
    u^o$ is well defined and in $\mathcal{D}^\eta$ for all $t \in
    \reali^+$ and $\eta > \eta_*$;
  \item $\tau_* (t,x) = \bar \tau_l$ and $v_* (t,x) = v_l (t)$ for all
    $(t,x) \in \reali^+ \times [0,m]$, with $v_l \in \W{1}{\infty}
    (\reali^+;\reali)$;
  \item for all $\bar x \in \reali$, the traces converge in the sense
    that
    \begin{displaymath}
      \lim_{\eta\to+\infty}
      \int_0^t \norma{u^\eta (s, \bar x \pm) - u_* (s, \bar x\pm)} \d{s}
      =
      0 \,;
    \end{displaymath}
  \item let $u_g = u_* \, \chi_{\strut\reali \setminus [0,m]}$. Then,
    the map $t \to \left(u_g (t), v_l (t)\right)$
    solves~\eqref{eq:IncompL} in the sense of
    Definition~\ref{def:solInc};
  \item as $\eta \to +\infty$, the map $p^{\eta} \circ \tau_l^{\eta}$
    weak* converges in $\L\infty$ to a function $p_l \in \L\infty
    (\reali^+ \times [0,m]; \reali^+)$ such that
    \begin{displaymath}
      \begin{array}{l}
        p_l (t,0+) \mbox{ is well defined and equals }
        p_g\left(\tau_g (t, 0-)\right)
        \mbox{ for a.e. }t \in \reali^+ \,,
        \\[4pt]
        p_l (t,m-) \mbox{ is well defined and equals }
        p_g\left(\tau_g (t, m+)\right)
        \mbox{ for a.e. }t \in \reali^+ \,,
        \\[4pt]
        p_l (t,x)
        =
        \left(1-\frac{x}{m}\right) \, p_g\left(\tau_g (t, 0-)\right)
        +
        \frac{x}{m} \, p_g\left(\tau_g (t, m+)\right)
        \mbox{ for a.e. }t \in \reali^+ \,.
      \end{array}
    \end{displaymath}
  \end{enumerate}
\end{theorem}

\noindent The proof of Theorem~\ref{thm:Main} is deferred to
Section~\ref{sec:TD}.

Observe that from the Eulerian coordinates' point of view, the
locations of the boundaries of the liquid phase can be recovered
through a time integration:
\begin{displaymath}
  \begin{array}{rcl@{\qquad\qquad}rcl}
    a^{\eta}(t) & = &
    \displaystyle
    a_{o} + \int_{0}^{t}v_{g}^{\eta}\left(s,0-\right)ds
    &
    a(t) & = &
    \displaystyle
    a_{o} + \int_{0}^{t}v_{l}\left(s\right)ds
    \\[5pt]
    b^{\eta}(t) & = &
    \displaystyle
    b_{o} + \int_{0}^{t}v_{g}^{\eta}\left(s,m+\right)ds
    &
    b(t) & = &
    \displaystyle
    b_{o} + \int_{0}^{t}v_{l}\left(s\right)ds \,.
  \end{array}
\end{displaymath}
The boundaries of the two phases turns out to be Lipschitz continuous
functions and moreover $a^{\eta}\to a$, $b^{\eta}\to b$ uniformly on
compact sets as $\eta\to+\infty$. Moreover, point~2.~in
Definition~\ref{def:solInc} justifies the usual relation between the
acceleration of the droplet and the pressure difference on its sides.

\section{Technical Details}
\label{sec:TD}

We collect below a few facts about he $p$-system in Lagrangian
coordinates
\begin{displaymath}
  \left\{
    \begin{array}{ll}
      \tau_t - v_x = 0
      \\
      v_t + \left[p(\tau)\right]_x = 0 \,.
    \end{array}
  \right.
  \quad t \in \left[0, +\infty \right[ \, , \; x \in \reali
\end{displaymath}
The eigenvalues are
\begin{equation}
  \label{eq:lambda}
  \lambda_1(\tau,v)
  =
  -\sqrt{-p'(\tau)}
  \qquad \qquad
  \lambda_2(\tau,v)
  =
  \sqrt{-p'(\tau)}
\end{equation}
so that the Lax shock and rarefaction curves in the $(v, \tau)$--plane
are
\begin{equation}
  \label{eq:Lax}
  \begin{array}{rcl@{\quad\mbox{if }}rcl}
    S_1(u, \sigma)
    & = &
    \left[
      \begin{array}{l}
        p-\sigma
        \\
        v - \sqrt{\left(\tau\left(p-\sigma\right)-\tau (p)\right)\sigma}
      \end{array}
    \right]
    &
    \sigma & < & 0
    \\[15pt]
    R_1 (u, \sigma)
    & = &
    \left[
      \begin{array}{l}
        p-\sigma
        \\
        \displaystyle
        v - \int_{p}^{p-\sigma} \sqrt{-\tau' (\pi)} \d\pi
      \end{array}
    \right]
    &
    \sigma & > & 0
    \\[18pt]
    S_2(u, \sigma)
    & = &
    \left[
      \begin{array}{l}
        p+\sigma
        \\
        v - \sqrt{-\left(\tau\left(p+\sigma\right)-\tau (p)\right)\sigma}
      \end{array}
    \right]
    &
    \sigma & < & 0
    \\[15pt]
    R_2 (u, \sigma)
    & = &
    \left[
      \begin{array}{l}
        p+\sigma
        \\
        \displaystyle
        v + \int_{p}^{p+\sigma} \sqrt{-\tau' (\pi)} \d\pi
      \end{array}
    \right]
    &
    \sigma & > & 0
  \end{array}
\end{equation}
Fix a state $\bar u_l$ and using the
pressure~\eqref{eq:linearPressure}, the Lax curves take the form
\begin{equation}
  \label{eq:linearLaxCurves}
  \mathcal{L}_1(u, \sigma)
  =
  \left[
    \begin{array}{l}
      p-\sigma
      \\
      v + \frac{\sigma}{\eta}
    \end{array}
  \right]
  \qquad \mbox{ and }\qquad
  \mathcal{L}_2(u, \sigma)
  =
  \left[
    \begin{array}{l}
      p+\sigma
      \\
      v + \frac{\sigma}{\eta}
    \end{array}
  \right] \,.
\end{equation}

\begin{proofof}{Theorem~\ref{thm:WFT}}
  In view of~\cite[Theorem~2.5]{ColomboSchleper}, it is sufficient to
  prove only the estimates at~\ref{it:thm1}.~and~\eqref{eq:new}. To
  this aim, we construct approximate solutions through an algorithm
  different from that in~\cite{ColomboSchleper}. Once a subsequence of
  these approximate solutions is proved to converge, by the properties
  of the Standard Riemann Semigroup~\cite{BressanLectureNotes}, it is
  known that the present approximations converge to the solution
  constructed in~\cite{ColomboSchleper}.

  \paragraph{Definition of the Algorithm.} Fix $\epsilon > 0$. We
  approximate the initial datum $u^o$ by a sequence $u^o_{\epsilon}$
  of piecewise constant initial data with a finite number of
  discontinuities such that $\norma{u^o_{\epsilon}-u^o}_{{\mathbf
      L}^1} \leq\epsilon$.  At each point of jump in the approximate
  initial condition, we solve the corresponding Riemann problem using
  the front tracking algorithm as stated
  in~\cite[Chapter~7]{BressanLectureNotes}.  We approximate each
  rarefaction wave by a rarefaction fan by means of (non-entropic)
  shock waves. Similarly to what happens in the general
  case~\cite[Chapter~5]{BressanLectureNotes}, there exists a constant
  $\delta^o > 0$ such that each of the above Riemann problems has a
  unique approximate solution as long as $\tv (u^o) < \delta^o$.

  This construction can be extended up to the first time $t_1$ at
  which two waves interact. At time $t_1$, the so constructed
  functions are piecewise constant with a finite number of
  discontinuities. We can thus iterate the previous construction at
  any subsequent interaction, provided suitable upper bounds on the
  total variation of the approximate solutions are available. These
  bounds essentially rely on the interaction estimates below.

  As it is usual in this contest, see~\cite{BressanLectureNotes}, we
  may assume that no more than 2 waves interact at any interaction
  point. Moreover, rarefaction waves, once arisen, are not further
  split even if their size exceeds the threshold $\epsilon$ after
  subsequent interactions, with other waves or with the phase
  boundaries. Besides, we call $\hat\lambda^\eta$ an upper bound for
  the moduli of the propagation speeds of all waves.

  Specific to the present construction, is our choice to parametrize
  the Lax curves as in~\eqref{eq:Lax} and, hence, waves' sizes are
  measured through the pressure variation $\sigma$ between the two
  states on the sides of the wave.

  \paragraph{Interaction Estimates.} We recall the classical Glimm
  interaction estimates, see~\cite[Chapter~7, formul\ae~(7.31)-
  (7.32)]{BressanLectureNotes}, which holds for any smooth
  parametrization:
  \begin{figure}[htpb]
    \centering
    \begin{psfrags}
      \psfrag{i+}{$\sigma_1^+$} \psfrag{i-}{$\sigma_1^-$}
      \psfrag{j+}{$\sigma_2^+$} \psfrag{j-}{$\sigma_2^-$}
      \psfrag{up}{$\sigma_1^+$} \psfrag{2+}{$\sigma_2^+$}
      \psfrag{ipp}{$\sigma_2''$} \psfrag{ip}{$\sigma_2'$}
      \includegraphics[width=8.5cm]{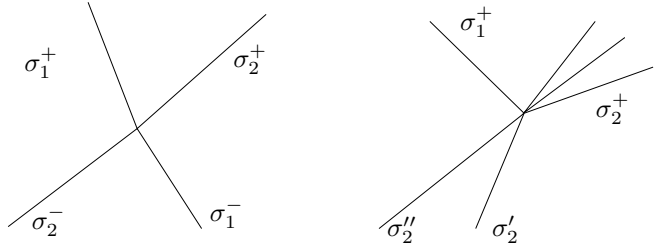}
    \end{psfrags}
    \caption{Left, an interaction between waves of different
      families. Right, an interaction between waves of the same
      family.\label{fig:interactions}}
  \end{figure}
  \begin{equation}
    \label{eq:ieStandard}
    \begin{array}{rcl}
      \modulo{\sigma_1^+ - \sigma_1^-}
      +
      \modulo{\sigma_2^+ - \sigma_2^-}
      & \leq &
      C \cdot \modulo{\sigma_1^- \sigma_2^-}
      \\
      \modulo{\sigma_1^+ - (\sigma'+\sigma'')}
      +
      \modulo{\sigma_2^+}
      & \leq &
      C \cdot \modulo{\sigma'\sigma''} 
    \end{array}
  \end{equation}
  where we used the notation described in
  Figure~\ref{fig:interactions}.  The estimates on the waves' sizes in
  the case of interactions involving the interface are as follows. In
  the case of Figure~\ref{fig:ie},
  \begin{equation}
    \label{eq:ieFig2}
    \modulo{\sigma_1^+} \leq C \, \modulo{\sigma_2^-}
    \qquad \mbox{ and } \qquad
    \modulo{\sigma_2^+} \leq C \, \modulo{\sigma_2^-}
  \end{equation}
  while in the case of Figure~\ref{fig:2},
  \begin{equation}
    \label{eq:ieFig3}
    \modulo{\sigma_1^+} + \modulo{\sigma_2^+} = \modulo{\sigma_2^-} \,.
  \end{equation}
  Note that the constant $C$ above depends only on quantities related
  to the gas phase, is uniformly bounded when $u_g$ varies in a
  compact set and, in particular, is independent from $\eta$.
  \begin{figure}[htpb]
    \centering
    \begin{psfrags}
      \psfrag{s1-}{$\sigma_2^-$} \psfrag{s1+}{$\sigma_1^+$}
      \psfrag{s2+}{$\sigma_2^+$} \psfrag{u1b}{$u_g$}
      \psfrag{u2b}{$u_l$} \psfrag{u1-}{$u_g^-$} \psfrag{u1+}{$u_g^+$}
      \psfrag{u2+}{$u_l^+$}
      \includegraphics[width=0.35\textwidth]{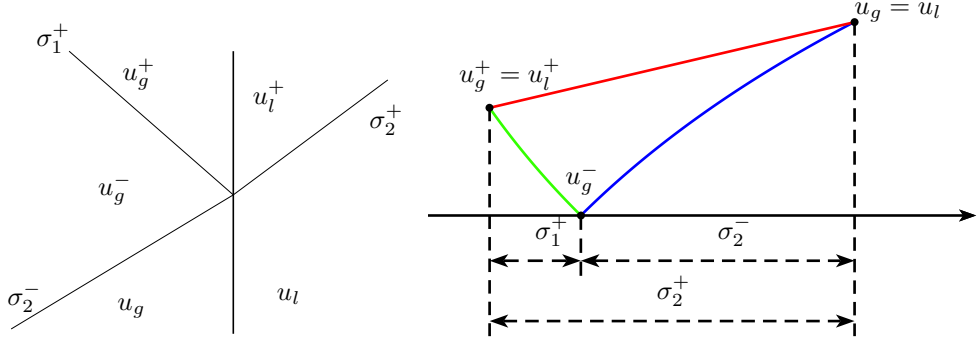}
    \end{psfrags}
    \newrgbcolor{ttffqq}{0.2 1 0} \newrgbcolor{xdxdff}{0.49 0.49 1}
    \psset{xunit=4.0cm,yunit=4.0cm,algebraic=true,dotstyle=o,
      dotsize=3pt 0,linewidth=0.8pt,arrowsize=3pt 2,arrowinset=0.25}
    \begin{pspicture*}(0.5,-0.4)(2.4,1.)
      \psaxes[labelFontSize=\scriptstyle,xAxis=true,yAxis=false,labels=none,Dx=10,Dy=10]{->}(0,0)(-0.1,-0.49)(2.3,0.8)
      \psplot[linewidth=1.pt,linecolor=ttffqq,plotpoints=200]{0.7}{1.}{-ln(x)}
      \psplot[linewidth=1.pt,linecolor=blue,plotpoints=200]{1.}{1.9}{ln(x)}
      \psline[linewidth=1.pt,linecolor=red](0.7,0.357)(1.9,0.641)
      \psline[linestyle=dashed,linewidth=1.pt](0.7,0.357)(0.7,-0.4)
      \psline[linestyle=dashed,linewidth=1.pt](1,0.)(1,-0.2)
      \psline[linestyle=dashed,linewidth=1.pt](1.9,0.641)(1.9,-0.4)
      \psline[linestyle=dashed,linewidth=1.pt]{<->}(1.,-0.15)(1.9,-0.15)
      \psline[linestyle=dashed,linewidth=1.pt]{<->}(0.7,-0.15)(1.,-0.15)
      \psline[linestyle=dashed,linewidth=1.pt]{<->}(0.7,-0.35)(1.9,-0.35)
      \psdots[dotstyle=*](1,0) \psdots[dotstyle=*](0.7,0.357)
      \psdots[dotstyle=*](1.9,0.641)
      \rput[bl](1.45,-0.1){$\sigma_{2}^{-}$}
      \rput[bl](0.85,-0.1){$\sigma_{1}^{+}$}
      \rput[bl](1.25,-0.3){$\sigma_{2}^{+}$}
      \rput[bl](0.95,0.07){$u_{g}^{-}$}
      \rput[bl](0.6,0.41){$u_{g}^{+}=u_{l}^{+}$}
      \rput[bl](1.9,0.641){$u_{g}=u_{l}$}
    \end{pspicture*}
    \\
    \caption{Notation used in the interaction estimates involving the
      phase interface.\label{fig:ie}}
  \end{figure}

  \begin{figure}[h]
    \centering \qquad
    \begin{psfrags}
      \psfrag{s1+}{$\sigma_2^+$} \psfrag{s1-}{$\sigma_1^-$}
      \psfrag{s2+}{$\sigma_1^+$} \psfrag{u1b}{$u_l$}
      \psfrag{u2b}{$u_g$} \psfrag{u1-}{$u_l^-$} \psfrag{u1+}{$u_l^+$}
      \psfrag{u2+}{$u_g^+$}
      \includegraphics[width=0.36\textwidth]{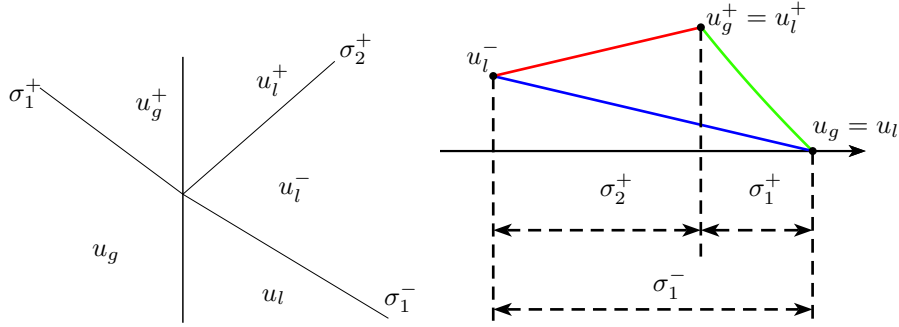}
    \end{psfrags}
    \newrgbcolor{ttffqq}{0.2 1 0} \newrgbcolor{xdxdff}{0.49 0.49 1}
    \psset{xunit=7.0cm,yunit=7.0cm,algebraic=true,dotstyle=o,
      dotsize=3pt 0,linewidth=0.8pt,arrowsize=3pt 2,arrowinset=0.25}
    \begin{pspicture*}(0.3,-0.33)(1.3,0.29) \psaxes[labelFontSize =
      \scriptstyle, xAxis = true, yAxis = false, labels = none, Dx =
      10, Dy = 10]{->}(0,0)(-0.1,-0.49)(1.1,0.8) \psplot[linewidth =
      1.pt, linecolor = ttffqq, plotpoints = 200]{0.79}{1.}{-ln(x)}
      \psline[linewidth=1.pt,linecolor=blue](1,0.)(0.4,0.142)
      \psline[linewidth=1.pt,linecolor=red](0.4,0.142)(0.79,0.2343)
      \psline[linestyle=dashed,linewidth=1.pt](0.4,0.142)(0.4,-0.4)
      \psline[linestyle=dashed,linewidth=1.pt](1,0.)(1,-0.4)
      \psline[linestyle=dashed,linewidth=1.pt](0.79,0.2343)(0.79,-0.2)
      \psline[linestyle=dashed,linewidth=1.pt]{<->}(0.4,-0.15)(0.79,-0.15)
      \psline[linestyle=dashed,linewidth=1.pt]{<->}(0.79,-0.15)(1.,-0.15)
      \psline[linestyle=dashed,linewidth=1.pt]{<->}(0.4,-0.3)(1.,-0.3)
      \psdots[dotstyle=*](1,0) \psdots[dotstyle=*](0.79,0.2343)
      \psdots[dotstyle=*](0.4,0.142)
      \rput[bl](0.6,-0.1){$\sigma_{2}^{+}$}
      \rput[bl](0.88,-0.1){$\sigma_{1}^{+}$}
      \rput[bl](0.7,-0.28){$\sigma_{1}^{-}$}
      \rput[bl](1,0.02){$u_{g}=u_{l}$}
      \rput[bl](0.80,0.22){$u_{g}^{+}=u_{l}^{+}$}
      \rput[bl](0.35,0.15){$u_{l}^{-}$}
    \end{pspicture*}
    \caption{Interaction estimate at the gas--liquid interface in the
      $(p,v)$ plane. Note that the two states $u_{g}^{+}$ and
      $u_{l}^{+}$ have different specific volumes but, due to the
      interface conditions, share the same pressure and the same
      speed.\label{fig:2}}
  \end{figure}

  \paragraph{Bounds on the Total Variation.} We follow the classical
  techniques based on Glimm functionals,
  see~\cite{BressanLectureNotes}. To this aim, introduce the
  quantities
  \begin{displaymath} \begin{array}{@{}r@{\,}c@{\,}l@{\qquad}r@{\,}c@{\,}l@{\qquad}r@{\,}c@{\,}l@{}}
      V_g^- & = & \displaystyle \sum_{i=1}^2 \sum_{x_\alpha<0}
      K_i^{g-} \modulo{\sigma_{i,\alpha}} & V_l & = & \displaystyle
      K_l \, \sum_{i=1}^2 \sum_{x_\alpha \in [0,m]}
      \modulo{\sigma_{i,\alpha}} & V_g^+ & = & \displaystyle
      \sum_{i=1}^2 \sum_{x_\alpha>m} K_i^{g+}
      \modulo{\sigma_{i,\alpha}}
      \\
      Q_g^- & = & \displaystyle \sum_{\mathcal{A}_g^-}
      \modulo{\sigma_{i,\alpha} \sigma_{j,\beta}} & & & & Q_g^+ & = &
      \displaystyle \sum_{\mathcal{A}_g^+} \modulo{\sigma_{i,\alpha}
        \sigma_{j,\beta}}
      \\
      \Upsilon_g^- & = & V_g^- + H \, Q_g^- & & & &\Upsilon_g^+ & = &
      V_g^+ + H \, Q_g^+
    \end{array}
  \end{displaymath}
  where $\mathcal{A}_g^-$, respectively $\mathcal{A}_g^+$, is the set
  of pairs of approaching waves both with support in $x<0$,
  respectively $x>m$. The weights $K_i^{g\pm}$, $K_l$ and $H$ are
  defined below. They are positive numbers independent from $\eta$.

  Define $\Upsilon=\Upsilon_{g}^{-}+V_{l}+\Upsilon_{g}^{+}$ and
  consider the various case:
  \begin{enumerate}[{Case}~1:]
    \setlength{\itemsep}{1pt} \setlength{\parskip}{0pt}
  \item An interaction in the interior of the liquid phase. then,
    thanks to the choice~\eqref{eq:linearPressure}, we have that
    $\Delta\Upsilon = \Delta V_l = 0$.
  \item \label{it:ie:2} An interaction at the interface. Consider
    first the case of Figure~\ref{fig:ie}. Then,
    \begin{eqnarray*}
      \Delta \Upsilon
      & = &
      \Delta \Upsilon_g^- + \Delta V_l
      \\
      & = &
      \Delta V_g^- + H \, \Delta Q_g^- + \Delta V_l
      \\
      & \leq &
      K_1^{g-} \, \modulo{\sigma_1^+}
      -
      K_2^{g-}\modulo{\sigma_2^-}
      +
      H \, \modulo{\sigma_1^+} \, \sum_{x_\alpha<0}\modulo{\sigma_\alpha}
      +
      K_l \, \modulo{\sigma_2^+}
      \\
      & \leq &
      \left(
        \left(
          K_1^{g-} + K_l + H \, \sum_{x_\alpha<0}\modulo{\sigma_\alpha}
        \right)C
        - K_2^{g-}
      \right)
      \modulo{\sigma_2^-} \,,
    \end{eqnarray*}
    while in the symmetric situation we get
    \begin{displaymath}
      \Delta \Upsilon \leq
      \left(
        \left(
          K_2^{g+} + K_l + H \, \sum_{x_\alpha>m}\modulo{\sigma_\alpha}
        \right)C
        - K_1^{g+}
      \right)
      \modulo{\sigma_1^-} \,.
    \end{displaymath}
    In the case of Figure~\ref{fig:2}:
    \begin{eqnarray*}
      \Delta \Upsilon
      & = &
      \Delta \Upsilon_g^- + \Delta V_l
      \\
      & = &
      \Delta V_g^- + H \, \Delta Q_g^- + \Delta V_l
      \\
      & \leq &
      K_1^{g-} \, \modulo{\sigma_1^+}
      +
      H \, \modulo{\sigma_1^+} \, \sum_{x_\alpha<0}\modulo{\sigma_\alpha}
      +
      K_{l} \, \modulo{\sigma_2^+}
      -
      K_{l} \, \modulo{\sigma_1^-}
      \\
      & = &
      K_1^{g-} \, \modulo{\sigma_1^+}
      +
      H \, \modulo{\sigma_1^+} \, \sum_{x_\alpha<0}\modulo{\sigma_\alpha}
      -
      K_{l} \, \modulo{\sigma_1^+}
      \\
      & \leq &
      \left(
        K_1^{g-} + H \, \sum_{x_\alpha<0}\modulo{\sigma_\alpha} - K_{l}
      \right) \modulo{\sigma_1^+} \,,
    \end{eqnarray*}
    while in the symmetric situation we get
    \begin{displaymath}
      \Delta \Upsilon
      \leq
      \left(
        K_2^{g+} + H \, \sum_{x_\alpha>m}\modulo{\sigma_\alpha} - K_{l}
      \right) \modulo{\sigma_2^+} \,,
    \end{displaymath}
  \item \label{it:ie:1} An interaction in the interior of the left gas
    phase. Then, the classical estimates~\eqref{eq:ieStandard} ensure
    that the standard Glimm functional decreases, i.e.,
    \begin{eqnarray*}
      \Delta\Upsilon
      & = &
      \Delta \Upsilon_g^-
      \\
      & = &
      \Delta V_g^- + H \Delta Q_g^-
      \\
      & \leq &
      (K_1^{g-} + K_2^{g-}) \, C  \, \modulo{\sigma_1^- \, \sigma_2^-}
      +
      H (C \delta -1) \, \modulo{\sigma_1^- \, \sigma_2^-}
      \\
      & \leq &
      \left(
        (K_1^{g-} + K_2^{g-}) \, C - H/2
      \right)
      \, \modulo{\sigma_1^- \, \sigma_2^-}
      \qquad \mbox{ and}
      \\
      \Delta\Upsilon
      & = &
      \Delta \Upsilon_g^-
      \\
      & = &
      \Delta V_g^- + H \Delta Q_g^-
      \\
      & \leq &
      (K_1^{g-} + K_2^{g-}) \, C  \, \modulo{\sigma' \, \sigma''}
      +
      H (C \delta -1) \,\modulo{\sigma' \, \sigma''}
      \\
      & \leq &
      \left(
        (K_1^{g-} + K_2^{g-}) \, C - H/2
      \right)
      \, \modulo{\sigma' \, \sigma''}
    \end{eqnarray*}
    in the two cases of Figure~\ref{fig:interactions}, provided
    $\delta$ is sufficiently small and $H$ is sufficiently
    large.  \end{enumerate}
  \noindent We now choose $K_1^{g-} = K_2^{g+} =1$, $K_l = 2$,
  $K_1^{g+} = K_2^{g-} = 4C$, $H = 4 (1+2C)C$ and finally $\delta =
  \min\left\{1/ (2C), 1/(2H)\right\}$, obtaining that
  $\Delta\Upsilon\leq 0$ at any interaction. Indeed, in the different
  cases considered above, we have:
  \begin{equation}
    \label{eq:DeltaUpsilon}
    \begin{array}{l@{\qquad}rcl}
      \mbox{Case~1:} & \Delta\Upsilon & = & 0
      \\
      \mbox{Case~2, Figure~\ref{fig:ie}, left boundary:}
      & \Delta\Upsilon & \leq & - \frac{C}{2} \, \modulo{\sigma_2^-}
      \\[5pt]
      \mbox{Case~2, right boundary:}
      & \Delta\Upsilon & \leq & - \frac{C}{2} \, \modulo{\sigma_1^-}
      \\
      \mbox{Case~2, Figure~\ref{fig:2}, left boundary:}
      & \Delta\Upsilon & \leq & -\frac{1}{2} \, \modulo{\sigma_1^+}
      \\[5pt]
      \mbox{Case~2, right boundary:}
      & \Delta\Upsilon & \leq & -\frac{1}{2} \, \modulo{\sigma_2^+}
      \\
      \mbox{Case~3, different families:}
      & \Delta\Upsilon & \leq & - C \, \modulo{\sigma_1^- \, \sigma_2^-}
      \\
      \mbox{Case~3, same family:}
      & \Delta\Upsilon & \leq & - C \modulo{\sigma' \, \sigma''} \,.
    \end{array}
  \end{equation}
  This implies that the map $t \to \Upsilon (t)$ decreases along the
  wave front tracking approximate solution, independently from
  $\epsilon$. This, thanks to the pressure
  law~\eqref{eq:linearPressure} and to the \emph{ad hoc} choice of the
  parametrization~\eqref{eq:linearLaxCurves}, leads to the following
  bounds uniform in $\epsilon$:
  \begin{equation}
    \label{eq:Uniform}
    \begin{array}{rclcl}
      \tv\left(p^\epsilon_l (t)\right)
      & \leq &
      \Upsilon (t)
      & \leq &
      \Upsilon (0)\,;
      \\
      \eta^2 \, \tv\left(\tau_l^\epsilon (t)\right)
      & \leq &
      \Upsilon (t)
      & \leq &
      \Upsilon (0) \,;
      \\
      \eta \, \tv\left(v_l^\epsilon (t)\right)
      & \leq &
      \Upsilon (t)
      & \leq &
      \Upsilon (0) \,.
    \end{array}
  \end{equation}

  \paragraph{Bounds on the Number of Interactions.}
  To be sure that the algorithm can be continued for all times one
  needs to prove that interaction times do not accumulate in finite
  time.

  Fix a positive $\epsilon$ and refer to the $\epsilon$--approximate
  wave front tracking solution $u^\epsilon$ defined above. Assume
  there exists a first time $t_\infty > 0$ such that the point
  $(t_\infty,x_\infty)$ is the limit of a sequence $(t_n, x_n)$ of
  interaction points, with $t_n < t_{n+1}$ for all $n$.

  \bigskip

  \noindent\textbf{Claim: $\boldsymbol{x_\infty \not\in \left]0, m\right[}$.}\quad By contradiction, assume that $x_\infty \in \left]0,
    m\right[$. Then, for suitable positive $\Delta t$ and $\Delta x$,
  the trapezoid
  \begin{displaymath}
    \mathcal{T} =
    \left\{
      (t,x) \in \reali^+ \times \reali \colon
      t \in [t_\infty-\Delta t, t_\infty]
      \mbox{ and }
      \modulo{x-x_\infty} \leq \Delta x + \hat\lambda^\eta (t_\infty-t)
    \right\}
  \end{displaymath}
  is contained in $\reali^+ \times [0,m]$. No waves can enter the
  sides of $\mathcal{T}$ and a finite number of waves is supported
  along the lower side of $\mathcal{T}$. Hence, only a finite number
  of interactions may take place inside $\mathcal{T}$ and no new waves
  are created therein, due to~\eqref{eq:linearPressure} in $\reali^+
  \times [0,m]$. This contradicts $(t_\infty, x_\infty)$ being the
  limit of an infinite sequence of distinct interaction points,
  proving the claim.

  We can thus assume that $x_\infty \geq m$.

\begin{figure}[h!]
  \centering
  \begin{psfrags}
    \psfrag{tinf}{$t_\infty$} \psfrag{tnxn}{$(t_n, x_n)$}
    \psfrag{inF}{$\mathcal{F}$} \psfrag{t0}{$t=0$}
    \includegraphics[width = 0.25\textwidth]{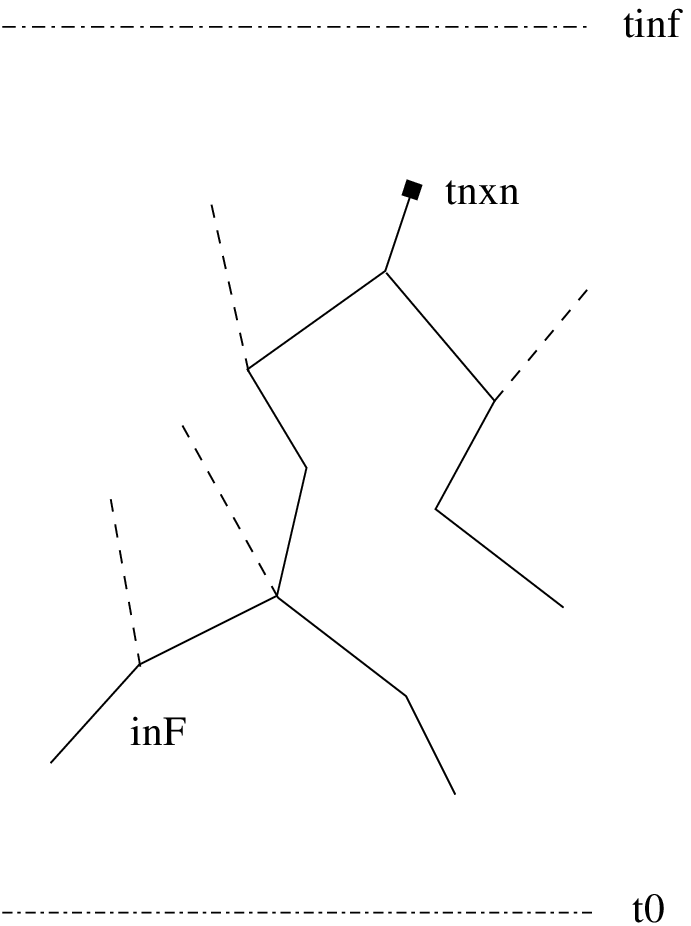}%
  \end{psfrags}
  \hfil
  \begin{psfrags}
    \psfrag{I0}{$\mathcal{I}_0$}
    \includegraphics[width = 0.25\textwidth]{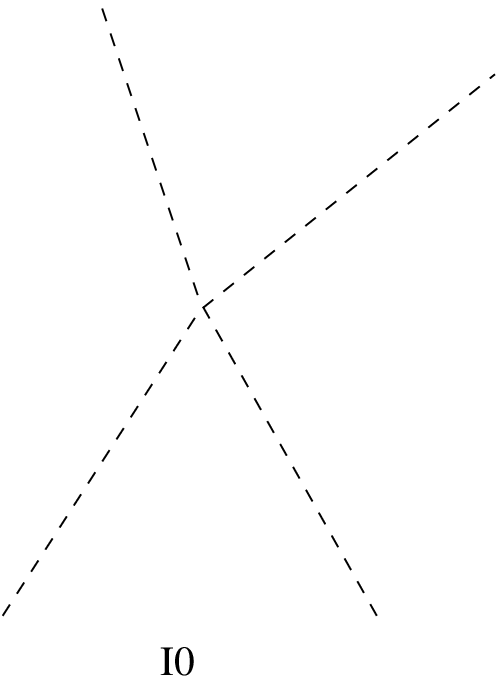}
  \end{psfrags}
  \caption{Left, an example of the set $\mathcal{F}$. Right, an
    interaction point in $\mathcal{I}_0$. Solid lines represent waves
    in $\mathcal{F}$, dashed ones are waves not in
    $\mathcal{F}$.\label{fig:F}}
\end{figure}
Call $\mathcal{F}$ the set of segments in $[0, t_\infty] \times
\left[m, +\infty\right[$ that support discontinuities in $u^\epsilon$
and can be connected forward in time along discontinuities in
$u^\epsilon$ to one of the points $(t_n, x_n)$, see
Figure~\ref{fig:F}.

Call $\mathcal{I}$ the set of all interaction points in $u^\epsilon$
in the strip $[0,t_\infty] \times \left[m, +\infty\right[$.
\begin{figure}[h!]
  \centering
  \begin{psfrags}
    \psfrag{tnxn}{$(t_n, x_n)$} \psfrag{I1}{$\mathcal{I}_1$}
    \includegraphics[width=0.30\textwidth]{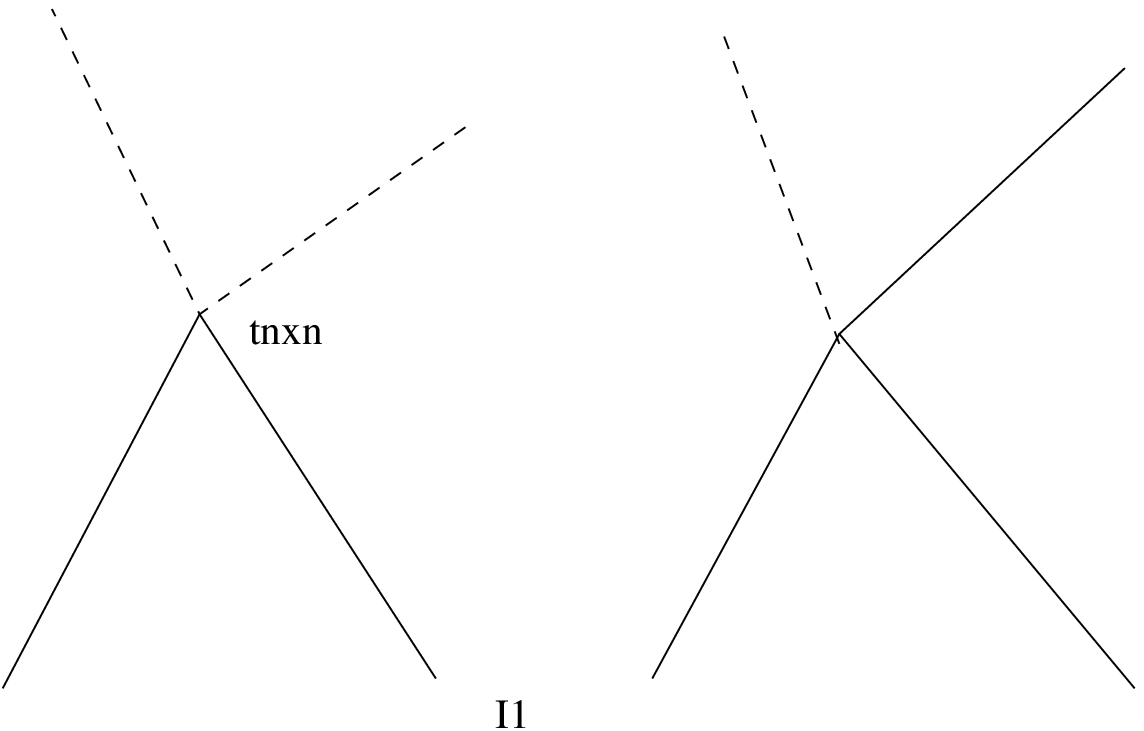}%
  \end{psfrags}
  \hfil\quad
  \begin{psfrags}
    \psfrag{I2}{$\mathcal{I}_2$} \psfrag{s1}{$\sigma''$}
    \psfrag{s2}{$\sigma'$}
    \includegraphics[width=0.16\textwidth]{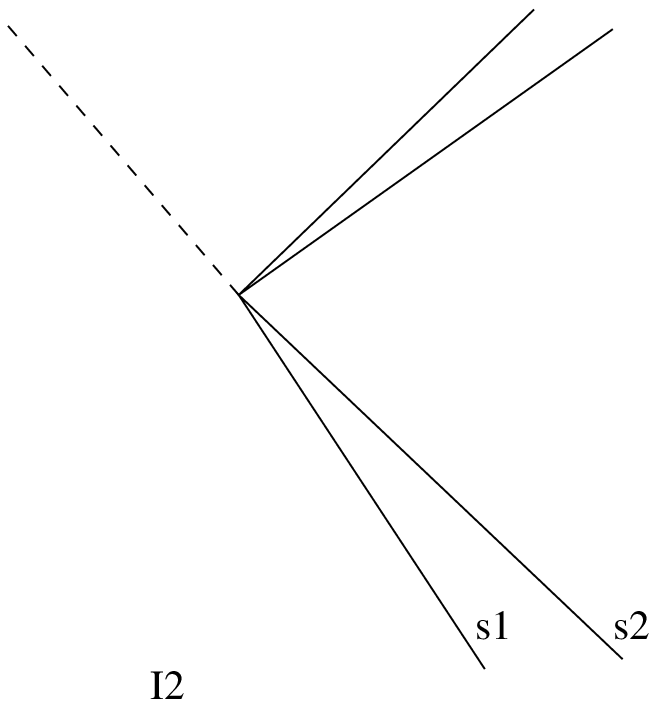}%
  \end{psfrags}
  \caption{Left and center, examples of interaction point in the set
    $\mathcal{I}_1$. Right, an interaction point in
    $\mathcal{I}_2$. Solid lines represent waves in $\mathcal{F}$,
    dashed ones are waves not in $\mathcal{F}$.\label{fig:I1I2}}
\end{figure}
With reference to figures~\ref{fig:I1I2} and~\ref{fig:I3I4}, remark
that $\mathcal{I}$ can be partitioned as follows:
\begin{eqnarray*}
  \mathcal{I}_0
  & = &
  \left\{
    (t,x) \in \mathcal{I} \colon
    x > m
    \mbox{ and no wave entering } (t,x) \mbox{ is in } \mathcal{F}
  \right\}
  \\
  \mathcal{I}_1
  & = &
  \left\{
    (t,x) \in \mathcal{I} \colon
    \begin{array}{@{}l}
      x > m
      \mbox{, at most }1 \mbox{ wave exiting } (t,x) \mbox{ is in }
      \mathcal{F}
      \\
      \mbox{and the two interacting waves are in } \mathcal{F}
    \end{array}
  \right\}
  \\
  \mathcal{I}_2
  & = &
  \left\{
    (t,x) \in \mathcal{I} \colon
    x > m
    \mbox{ and at least }2 \mbox{ outgoing waves of the same family are in }
    \mathcal{F}
  \right\}
  \\
  \mathcal{I}_3
  & = &
  \left\{
    (t,x) \in \mathcal{I} \colon
    \begin{array}{@{}l}
      x = m \mbox{, the wave entering }(t,x) \mbox{ from }
      \reali^+ \times \left[m, +\infty\right[ \mbox{ is in } \mathcal{F}
      \\
      \mbox{and the outgoing wave is not in } \mathcal{F}
    \end{array}
  \right\}
  \\
  \mathcal{I}_\infty
  & = &
  \left\{
    (t,x) \in \mathcal{I} \colon
    \begin{array}{@{}l}
      x > m
      \mbox{, the 2 waves exiting } (t,x)
      \mbox{ are one of the first family,}
      \\
      \mbox{one of the second and both belong to }\mathcal{F}
    \end{array}
  \right\}
  \\
  & &
  \cup
  \left\{
    (t,x) \in \mathcal{I} \colon
    x = m \mbox{ and the outgoing wave in }
    \reali^+ \times \left[m, +\infty\right[ \mbox{ is in } \mathcal{F}
  \right\}
\end{eqnarray*}
\begin{figure}[h!]
  \centering
  \begin{psfrags}
    \psfrag{I4}{$\mathcal{I}_3$} \psfrag{xm}{\hspace{-3mm}$x=m$}
    \includegraphics[width=0.12\textwidth]{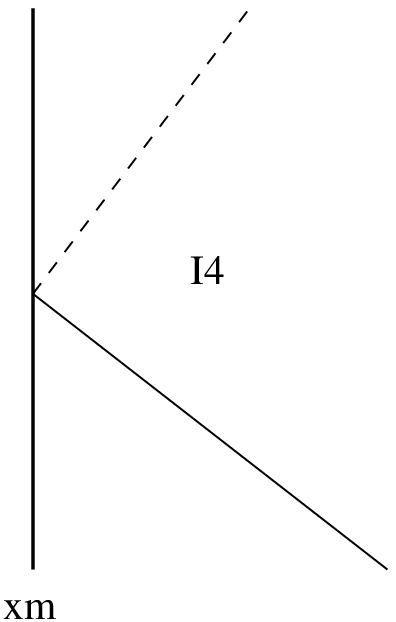}%
  \end{psfrags}
  \hfil
  \begin{psfrags}
    \psfrag{Iinf}{$\mathcal{I}_\infty$}
    \includegraphics[width=0.176\textwidth]{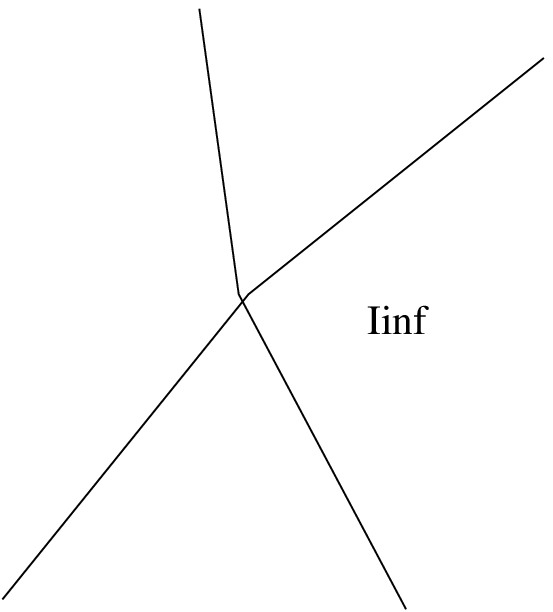}%
  \end{psfrags}
  \hfil
  \begin{psfrags}
    \psfrag{I3}{$\mathcal{I}_\infty$} \psfrag{xm}{\hspace{-3mm}$x=m$}
    \includegraphics[width=0.12\textwidth]{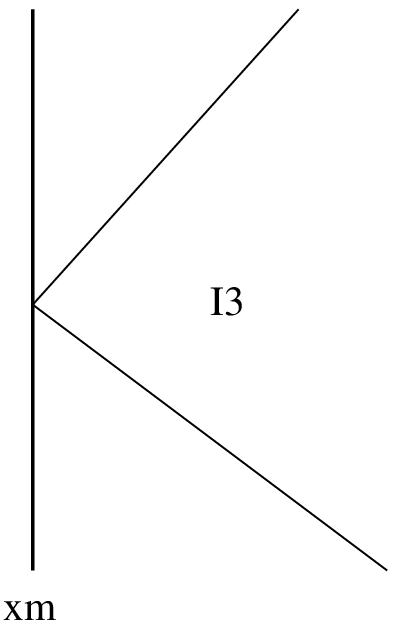}%
  \end{psfrags}
  \caption{Left, an example of interaction point in
    $\mathcal{I}_3$. Center, respectively right, an interaction point
    in $\mathcal{I}_\infty$ at $x>m$, respectively along the
    boundary. Solid lines represent waves in $\mathcal{F}$, dashed
    ones are waves not in $\mathcal{F}$, the thick line is the (right)
    phase boundary.\label{fig:I3I4}}
\end{figure}

\medskip

  \noindent\textbf{Claim: The set $\mathcal{I}_2$ is finite.}\quad Refer to Figure~\ref{fig:I1I2}, right. When two waves of the same family arise, they are rarefactions with total size bigger than $\epsilon$. Hence, the interacting waves are
  shocks of the same family with sizes $\sigma'$, $\sigma''$
  satisfying $\modulo{\sigma'\,\sigma''} > \epsilon/C$. Hence,
  $\Delta\Upsilon < - \epsilon$ and this can happen only a finite
  number of times.

  \medskip

  \noindent\textbf{Claim: The sets $\mathcal{I}_1$ and $\mathcal{I}_3$
    are finite.}\quad Call $\mathcal{Z} (\tau)$ the number of segments
  in $\mathcal{F}$ that intersect the line $t=\tau$. $\mathcal{Z}$ may
  increase only at the interactions in $\mathcal{I}_2$, hence a finite
  number of times. Then, at each interaction in $\mathcal{I}_1$ and in
  $\mathcal{I}_3$ it decreases at least by $1$. Thus, $\mathcal{I}_1$
  and $\mathcal{I}_3$ are finite.

  \medskip

  Thus, in the partition above of $\mathcal{I}$, only the sets
  $\mathcal{I}_0$ and $\mathcal{I}_\infty$ can be infinite. Let
  $\tau_\infty \in \left]0, t_\infty\right[$ be such that all the
  interactions in $\mathcal{I}$ in the time interval
  $\left[\tau_\infty, t_\infty\right[$ are in $\mathcal{I}_0$ or in
  $\mathcal{I}_\infty$.

  \medskip
  \noindent\textbf{Claim: $\boldsymbol{x_\infty>m}$.}\quad Let $t_n
  \in \left]\tau_\infty, t_\infty\right[$. Then, $x_n \geq m$ and it
  is possible to follow waves in $\mathcal{F}$ only of the second
  family converging to $x_\infty$. Waves of the second family have
  positive speed, hence $x_\infty > x_n \geq m$.

  \medskip
  \noindent\textbf{Claim: $\boldsymbol{(t_\infty,x_\infty)}$ does not
    exist.} Let $(t_n,x_n)$ be sufficiently near to
  $t_\infty,x_\infty$, so that $x_n-m > \hat\lambda^\eta
  (t_\infty-t_n)$ and $t_n > \tau_\infty$. Then, the phase boundary
  $x=m$ can not be reached following waves of the first family in
  $\mathcal{F}$ that exit $(t_n,x_n)$. On the other side, by the
  definition of $\mathcal{F}$, $(t_\infty,x_\infty)$ should be reached
  starting from $(t_n, x_n)$ following waves only of the first and
  only of the second family in $\mathcal{F}$ that connect points in
  $\mathcal{I}_\infty$. However, following waves of the first family
  one reaches a strictly decreasing sequence of points $x_{n_j}$, so
  that $x_\infty = \inf_j x_{n_j} < x_n$, whereas along waves of the
  second family one reaches a strictly increasing sequence of points
  $x_{n_i}$, so that $x_\infty = \sup_i x_{n_i} >x_n$. This
  contradicts the existence of $(t_\infty, x_\infty)$.

  \bigskip

  \paragraph{Lipschitz Continuity in Time.} Let $u^\epsilon (t) =
  (p^\epsilon, v^\epsilon) (t)$ denote the $\epsilon$-solution
  constructed above.  Introduce the following upper bound for the
  characteristic speeds in the gas phase $\Lambda =
  \norma{\sqrt{p'}}_{\C0}$ and note that the characteristic speed in
  the liquid phase is $\eta$.  Then, by the above definition of the
  approximate solution and by the
  parametrization~\eqref{eq:linearLaxCurves}, see
  also~\cite[Chapter~7, formula~(7.9)]{BressanLectureNotes}, if $t_1 <
  t_2$,
  \begin{eqnarray}
    \nonumber
    \norma{v^\epsilon (t_1) - v^\epsilon (t_2)}_{\L1}
    & \leq &
    \left(
      \Upsilon_g^-\left(u^\epsilon (t_1)\right)
      +
      \Upsilon_g^+\left(u^\epsilon (t_1)\right)
    \right) \, \Lambda \, \modulo{t_1-t_2}
    +
    V_l \left(u^\epsilon (t)\right) \, \modulo{t_1-t_2}
    \\
    \label{eq:LipTv}
    & \leq &
    \Upsilon (u^\epsilon_o) (1+\Lambda) \modulo{t_1-t_2} \,.
  \end{eqnarray}
  Passing now to the pressure, the same computations lead to:
  \begin{eqnarray*}
    \norma{p^\epsilon (t_1) - p^\epsilon (t_2)}_{\L1 (\reali\setminus[0,m])}
    & \leq &
    \left(
      \Upsilon_g^-\left(u^\epsilon (t_1)\right)
      +
      \Upsilon_g^+\left(u^\epsilon (t_1)\right)
    \right) \, \Lambda \, \modulo{t_1-t_2}
    \\
    \norma{p^\epsilon (t_1) - p^\epsilon (t_2)}_{\L1 ([0,m])}
    & \leq &
    \eta \, V_l \left(u^\epsilon (t)\right) \, \modulo{t_1-t_2}
  \end{eqnarray*}
  and the parametrization~\eqref{eq:Lax}--\eqref{eq:linearLaxCurves}
  give the following bounds on the specific volume
  \begin{eqnarray}
    \nonumber
    \norma{\tau^\epsilon_g (t_1) - \tau^\epsilon_g (t_2)}_{\L1 (\reali\setminus[0,m])}
    & \leq &
    C
    \left(
      \Upsilon_g^-\left(u^\epsilon (t_1)\right)
      +
      \Upsilon_g^+\left(u^\epsilon (t_1)\right)
    \right) \, \Lambda \, \modulo{t_1-t_2}
    \\
    \label{eq:LipTtau1}
    & \leq &
    C \, \Upsilon (u^\epsilon_o) \Lambda \modulo{t_1-t_2}
    \\
    \label{eq:LipTtau2}
    \norma{\tau^\epsilon_l (t_1) - \tau^\epsilon_l (t_2)}_{\L1 ([0,m])}
    & \leq &
    \frac{1}{\eta}
    V_l \left(u^\epsilon (t)\right) \, \modulo{t_1-t_2}
    \; \leq \;
    \frac{1}{\eta} \, \Upsilon (u^\epsilon_o) \modulo{t_1-t_2} \,.
  \end{eqnarray}

  \paragraph{The Limit.} Let $\epsilon_k$ be a sequence converging to
  $0$ as $k \to +\infty$. Then, the sequence $u^{\epsilon_k}$ of
  $\epsilon_k$-wave front tracking solutions satisfies Helly
  Compactness Theorem. By the uniqueness of the Riemann Semigroup
  constructed in~\cite[Theorem~2.5]{ColomboSchleper}, any of its
  converging subsequences converge to the unique semigroup of entropy
  solutions constructed therein. Then, the $\L1$-lower semicontinuity
  of the total variation allows to pass the bounds~\eqref{eq:Uniform}
  to the limit $k\to +\infty$, proving both inequalities
  in~\eqref{eq:new}.

  Concerning point~1, we have to show that there is a constant
  $\delta_{g}$ independent from $\eta$ and a positive
  $\delta_{l}^{\eta}$, such that all functions in the domain
  \begin{displaymath}
    \left\{u \in \bar
      u + (\L1\cap\BV) (\reali; \reali^+ \times \reali)\colon \tv
      (u_g) < \delta_g \mbox{ and } \tv (u_l) < \delta_l^\eta
    \right\}
  \end{displaymath}
  satisfy $\Upsilon < \delta$. By standard properties of the Riemann
  problem, it is easy to show that if the total variation of the
  initial data is sufficiently small, independently from $\eta$, then
  $\Upsilon_{g}^{-} + \Upsilon_{g}^{+} < \delta /2$.  Concerning the
  liquid phase, the choice~\eqref{eq:linearPressure} allows to compute
  the exact solution to any Riemann problem and obtain the estimate
  \begin{displaymath}
    V_{l}
    =
    K_l \, \tv \left(p^{\eta}\left(\tau_{l}(0+)\right)\right)
    =
    \frac{K_l}{2}
    \left[
      \eta^{2}\tv\left(\tau_{l}(0)\right) + \eta\tv\left(v_{l}(0)\right)
    \right]
  \end{displaymath}
  where we have also used \eqref{eq:linearPressure}.  Therefore,
  choosing $\delta^\eta_l < \delta/ (K_l \, \eta^2)$ implies
  $\Upsilon<\delta$.

  \paragraph{Point~\ref{it:thm:last}.}  By~\textbf{(p)}
  and~\eqref{eq:lambda}, the line $x=\bar x$ is noncharacteristic. Fix
  $\eta$, $\epsilon$ and a point $x > m$, the case $x<0$ being
  entirely similar. The estimate we now prove is an analogous
  to~\cite[Formula~(14.5.19)]{DafermosBook} in the present setting.

  Introduce the functional
  \begin{displaymath}
    \Xi_x (t)
    =
    \tv\left(p_g\left(\tau_g^{\epsilon,\eta} (\cdot, x)_{\strut\big|[0,t]}\right)\right)
    +
    \sum_{x_\alpha \in \left]m,x\right[} \modulo{\sigma_{2,\alpha}}
    +
    \sum_{x_\alpha \in \left]m, +\infty\right[} \modulo{\sigma_{1,\alpha}}
    +
    4 \, \Upsilon (t)
  \end{displaymath}
  and observe that it is non increasing in time. Indeed, considering
  all the possible interactions as above, we have:
  \begin{enumerate}[{Case}~1:]
  \item Clearly, $\Delta \Xi_x =0$.
  \item When the interaction is against the left boundary, in both
    cases of Figure~\ref{fig:ie} and Figure~\ref{fig:2}, clearly
    $\Delta \Xi_x = \Delta \Upsilon < 0$. When the interaction is as
    in against the right boundary and the incoming wave comes from the
    gas phase, by~\eqref{eq:ieFig2} and~\eqref{eq:DeltaUpsilon}, we
    have
    \begin{displaymath}
      \Delta\Xi
      =
      \modulo{\sigma_2^+} + 4\Delta \Upsilon
      \leq
      -C \modulo{\sigma_1^-}
      < 0 \,.
    \end{displaymath}
    When the interaction is against the right boundary and the
    incoming wave comes from the liquid phase, by~\eqref{eq:ieFig3}
    and~\eqref{eq:DeltaUpsilon}, we have:
    \begin{displaymath}
      \Delta\Xi
      =
      \modulo{\sigma_2^+} + 4\,\Delta \Upsilon
      \leq
      - \modulo{\sigma_2^+}
      < 0 \,.
    \end{displaymath}
  \item If the interaction is in the interior of the left gas phase,
    then clearly $\Delta \Xi_x = \Delta \Upsilon < 0$. In the case of
    an interaction between waves of different families at an
    interaction point $\bar x \in \left]m, x\right[$,
    by~\eqref{eq:ieStandard} and~\eqref{eq:DeltaUpsilon}, we have:
    \begin{displaymath}
      \Delta\Xi
      =
      \modulo{\sigma_2^+}
      -
      \modulo{\sigma_2^-}
      -
      4 \, C \, \modulo{\sigma_1^- \, \sigma_2^-}
      \leq
      -3\, C \, \modulo{\sigma_1^- \, \sigma_2^-}
      <0 \,.
    \end{displaymath}
    If $\bar x >x$, then, by~\eqref{eq:ieStandard}
    and~\eqref{eq:DeltaUpsilon},
    \begin{displaymath}
      \Delta\Xi
      =
      \modulo{\sigma_1^+}
      -
      \modulo{\sigma_1^-}
      -
      4 \, C \, \modulo{\sigma_1^- \, \sigma_2^-}
      \leq
      -3\, C \, \modulo{\sigma_1^- \, \sigma_2^-}
      <0
    \end{displaymath}
    and entirely analogous estimates hold when the interacting waves
    belong to the same family.
  \item We have now to consider also the case of a wave passing
    through $x$. Then, by the definition of $\Xi$, $\Delta \Xi = 0$.
  \end{enumerate}

  \noindent The monotonicity of $\Xi$ ensures the following estimate
  on the total variation of the $\epsilon$-wave front tracking
  approximate solution constructed above:
  \begin{displaymath}
    \tv\left(u_g^{\epsilon, \eta} (\cdot, x)\right)
    \leq
    \mathcal{K}_0 \, \Xi (0)
    \leq
    \mathcal{K}_1 \, \Upsilon \left(u^{\epsilon, \eta} (0+)\right)
    \leq
    \mathcal{K}_2 \, \Delta
  \end{displaymath}
  where $\mathcal{K}_1, \mathcal{K}_2$ are constants independent from
  the initial data, from $\eta$ and from $\epsilon$. Fix now $x_1$ and
  $x_2$ in the same gas phase. Then, similarly
  to~\cite[Formula~(14.4.7)]{DafermosBook},
  \begin{eqnarray*}
    \int_0^t \norma{u_g^{\epsilon, \eta} (s, x_2) - u_g^{\epsilon, \eta} (s,x_1)} \d{s}
    & \leq &
    \frac{1}{\inf \sqrt{-p_g'}}
    \left(
      \sup_{x \in \reali\setminus[0,m]} \tv\left(u_g^{\epsilon, \eta} (\cdot, x)\right)
    \right) \, \modulo{x_2-x_1}
    \\
    & \leq &
    \mathcal{K}_3 \, \Delta \, \modulo{x_2-x_1}
  \end{eqnarray*}
  for a suitable $\mathcal{K}_3$ independent from the initial data,
  from $\eta$ and from $\epsilon$. By the same arguments used in the
  paragraph above, we may pass to the limit $\epsilon \to 0$ obtaining
  for a.e.~$x_1,x_2$ the Lipschitz estimate \eqref{eq:uno}, completing
  the proof of Theorem~\ref{thm:WFT}.
\end{proofof}

\begin{proofof}{Theorem~\ref{thm:Main}}
  Point~1.~is a direct consequence of~\ref{it:thm1} in
  Theorem~\ref{thm:WFT}, with $\delta = \delta_g$.

  Let $\eta_k$ be a sequence with $\lim_{k\to+\infty} \eta_k =
  +\infty$. Helly Compactness Theorem~\cite[Chapter~2,
  Theorem~2.4]{BressanLectureNotes} can be applied thanks to the
  Lipschitz estimates~\ref{it:thm:LipT}.~in Theorem~\ref{thm:WFT} and
  ensures that a subsequence of $S^{\eta_k}_t u$ converges a.e.~to a
  function $u_* \in \C{0,1} \left(\reali^+;\bar u + \L1 (\reali;
    \reali^+ \times \reali)\right)$, where $u_* = (\tau_*, v_*)$,
  proving point~2.

  By~\eqref{eq:new}, for all $t \in \reali^+$, the function $u_*$ is
  constant for $x \in [0,m]$. Hence, the bounds
  at~\ref{it:thm:LipT}.~in Theorem~\ref{thm:WFT} ensure that $\tau^*
  (t,x) = \bar\tau$ for all $(t,x) \in \reali^+ \times [0, m]$ and
  that $v \in \W{1}{\infty} (\reali^+; \reali)$, completing the proof
  of point~3.

  To prove point~4., the case $\bar x \in [0,m]$ is immediate by the
  bounds~\eqref{eq:new} which ensure that the limit is independent
  from $x$ in the liquid phase. Assume that $\bar x \in \reali
  \setminus [0,m]$ and note that passing to the limit $\eta \to
  +\infty$ in~\ref{it:thm:last} of Theorem~\ref{thm:WFT}, we have
  \begin{equation}
    \label{eq:due}
    \int_0^t \norma{u_* (s, x_2) - u_* (s,x_1)} \d{s}
    \leq
    L \, \modulo{x_2-x_1}
  \end{equation}
  for a.e.~$x_1,x_2$ in the same gas phase. Consider the case of the
  right trace, the other case being entirely similar.  Let $x_n$ be a
  sequence converging to $\bar x$ from the right of points
  where~\ref{it:thm:last}.~in Theorem~\ref{thm:WFT}
  applies. Using~\eqref{eq:due},
  \begin{eqnarray*}
    & &
    \int_0^t \norma{u_g^\eta (s, \bar x +) - u_* (s, \bar x+)} \d{s}
    \\
    & \leq &
    \int_0^t \norma{u_g^\eta (s, \bar x +) - u_g^\eta (s, x_n)} \d{s}
    +
    \int_0^t \norma{u_g^\eta (s, x_n) - u_* (s, x_n)} \d{s}
    \\
    & &
    +
    \int_0^t \norma{u_* (s, x_n) - u_* (s, \bar x+)} \d{s}
    \\
    & \leq &
    2\, L \modulo{\bar x - x_n}
    +
    \int_0^t \norma{u_g^\eta (s, x_n) - u_* (s, x_n)} \d{s}
  \end{eqnarray*}
  To complete the proof of point~4., pass now to the $\limsup$ as
  $\eta \to +\infty$ and then to the limit as $n \to +\infty$.

  To prove point~5., note that the regularity
  condition~\eqref{eq:regularity} clearly holds. Point~1.~in
  Definition~\ref{def:solInc} is proved passing to the limit in the
  definition of weak entropy solution, which is possible by the
  convergence proved above.  Recall now that the integral formulation
  of~\eqref{eq:ModelL} implies
  \begin{displaymath}
    \int_0^m \left(v_l^\eta (t,x) - \bar v_l\right) \d{x}
    =
    \int_0^t \left(
      p_g\left(\tau_g^\eta (s, 0-)\right)
      -
      p_g\left(\tau_g^\eta (s, m+)\right)
    \right)
    \d{s}
  \end{displaymath}
  and, thanks to the convergence of the traces proved above, in the
  limit $\eta\to +\infty$ we have
  \begin{displaymath}
    v_l (t) - \bar v_l
    =
    \int_0^t \frac{1}{m}
    \left(
      p_g \left(\tau_g(s, 0-)\right) - p_g \left(\tau_g(s, m+)\right)
    \right) \d{s} \,.
  \end{displaymath}
  which is the integral formulation of point~2.~in
  Definition~\ref{def:solInc}. Finally, point~3.~in the same
  definition immediately follows from the convergence of the traces
  proved above.

  To prove point~6., define
  \begin{displaymath}
    \pi^\eta (t,x)
    =
    \left\{
      \begin{array}{l@{\quad}rcl}
        p_g\left(\tau_g^\eta (t,x)\right)
        & x & \in & [0,m]
        \\[6pt]
        p^\eta\left(\tau_{l\vphantom{g}}^\eta (t,x)\right)
        & x & \in & \reali \setminus [0,m]
      \end{array}
    \right.
  \end{displaymath}
  and let $\eta_n$ be an arbitrary real sequence converging to
  $+\infty$. Note that the sequence $\pi^{\eta_n}$ is uniformly
  bounded in $\L\infty$ by~\eqref{eq:new}, hence it admits a
  subsequence which is weak* convergent to a limit $\pi \in \L\infty
  (\reali^+ \times \reali; \reali^+$). Passing to the limit $n \to
  +\infty$ in the definition of weak solution, we have
  \begin{equation}
    \label{eq:integral}
    \int_{\reali^+} \int_{\reali}
    \left(
      v_* (t,x) \, \partial_t \phi (t,x) + \pi (t,x) \, \partial_x \phi (t,x)
    \right)
    \d{x} \d{t} = 0
  \end{equation}
  for all $\phi \in \Cc1 (\reali^+ \times \reali; \reali)$. Call now
  $p_l = \pi_{|[0,m]}$ and remark that in $[0,m]$ the following
  equality is satisfied in the sense of distributions:
  \begin{displaymath}
    \partial_t v_l + \partial_x p_l = 0
  \end{displaymath}
  showing that, by~3., $\partial_x p_l(t, \cdot)$ is constant in $x$
  and, hence, the map $x \to p_l (t,x)$ is linear. The existence of
  the traces immediately follows. Moreover, by~\eqref{eq:integral},
  necessarily $\pi (t,0-) = p_l (t, 0+) $ and $\pi (t,m+) = p_l (t,
  m-)$, which shows also the linear interpolation formula and, hence,
  that $p_l$ is independent from the choice of the sequence $\eta_n$
\end{proofof}

\medskip

\noindent\textbf{Acknowledgment:} The authors acknowledge the support by the German Research Foundation (DFG) in the framework of the Collaborative
Research Center Transregio 75 \emph{Droplet Dynamics under Extreme
  Ambient Conditions}, by the Elite program for postDocs of the
Baden-W\"urttemberg Foundation, by the PRIN-2009 project \emph{Fluid
  dynamics hyperbolic Equations and conservation Laws} and by the
research program \emph{Conservation Laws and Applications} of the
University of Brescia.

{\small

  \bibliographystyle{abbrv}

  \bibliography{c2i}

}

\end{document}